\documentclass[12pt,A4]{article}
\usepackage{amssymb}
\usepackage{amsmath}

\newtheorem{theorem}{Theorem}

\newtheorem{lemma}[theorem]{Lemma}

\newtheorem{proposition}[theorem]{Proposition}
\newtheorem{remark}[theorem]{Remark}

\def\func#1{\mathop{\rm #1}}

\begin{document}

\title{
SOME CHARACTERIZATIONS OF THE SPHERICAL HARMONICS COEFFICIENTS FOR
ISOTROPIC RANDOM FIELDS}
\author{Paolo Baldi\\
Dipartimento di Matematica\\
Universit\`a di Roma - Tor Vergata\\
Via della Ricerca Scientifica\\
00133 Roma, Italy.\\
baldi@mat.uniroma2.it
\and
Domenico Marinucci\\
Dipartimento di Matematica\\ Universit\`a di Roma - Tor Vergata\\
Via della Ricerca Scientifica\\ 00133 Roma, Italy.\\
marinucc@mat.uniroma2.it
}
\date{}
\maketitle
\begin{abstract}
\noindent In this paper we provide some simple characterizations for the spherical harmonics coefficients
of an isotropic random field on the sphere. The main result is a characterization of
isotropic gaussian fields through independence of the coefficients of their development in
spherical harmonics.
\end{abstract}
\smallskip


\section{Introduction}
In recent years the study of real random fields on the sphere has received much
attention,
as this topic is a necessary tool in the statistical study of the CMB (Cosmic Microwave Background)
radiation.
The existing physical literature is huge, with particular emphasis on testing for
Gaussianity and isotropy (e.g., for a small sample, \cite{DCB:2003}, \cite{cruz:05}, \cite{Pa:2004})
as both these issues have deep implications for cosmological physics.
See for reviews on this subject \cite{MR2106384} and \cite{MA:2005}, where many more references
can be found. We also mention \cite{JSDAF:05}, \cite{AGMNW:04} and \cite{MR2065205} for a mathematical treatment of the subject.

A natural tool for this kind of enquiry is the development of the random field in a
series of spherical harmonics. This raises some simple, but not so obvious questions.
For an isotropic random field $T$, given its development in spherical harmonics (see \S 2),
\begin{equation}\label{develop}
T(x)=\sum_{\ell=0}^\infty\sum_{m=-\ell}^\ell a_{\ell m} Y_{\ell m}(x)
\end{equation}
what properties must be expected to be satisfied by the coefficient $a_{\ell m}$? Or,
from a different point of view, which conditions must be verified by the coefficients
in order that the development above defines an isotropic random field?

These questions are of interest both from a statistical point of
view (for instance in order to devise a test suitable to detect
non-Gaussianity or anisotropy) or from a probabilistic point of
view (how to sample an isotropic random field).

In this paper we provide some results in this perspective. More
precisely, first we prove that, for an isotropic random field, the
coefficients are necessarily uncorrelated. This fact is well
known, but we give a proof that holds without the assumption that
the field is mean square continuous.

Then it is proved that each of the complex r.v.'s $a_{\ell m}$ has
a distribution whose phase is uniform on $[-\pi,\pi]$. This
implies in particular that, for any isotropic field, the ratio
$\func{Re} a_{\ell m}/\func{Im} a_{\ell m}$ is necessarily
distributed accordingly to a Cauchy distribution.

Finally, having remarked that if the field $T$ is gaussian then
the $a_{\ell m}$'s, $\ell=0,1,\dots$, $m=0,\dots,\ell$ are
independent, we prove that also the converse is true. Thus the
only isotropic fields such that the $a_{\ell m}$'s,
$\ell=0,1,\dots$, $m=0,\dots,\ell$ are independent are those that
are gaussian. This, which is the main result of this paper, is not
a consequence of the central limit theorem, but results from some
classical characterization of gaussian random variables, through
independence of some linear statistics.

This result gives a rigorous proof of claims that can be occasionally found
in the cosmological literature (see \cite{CoMa:2001} e.g.)

\section{Isotropic random fields}
In this section we recall some well known facts about isotropic random
fields $T$ defined on the unit sphere $S^{2}=\{ x\in
{\mathbb R}^{3}:| x | =1\}$. We assume that these fields
are isotropic in the strong sense,
that is, their probability law is invariant with respect to the action of $
SO(3)$. The isotropy assumption can then be stated as follows: for all $g\in SO(3)$
and $x_{1},\dots ,x_{p}\in S^{2}$, the two vectors %
\begin{equation*}
(T(gx_{1}),\dots ,T(gx_{p}))\quad \mbox{and}\quad(T(x_{1}),\dots ,T(x_{p}))\text{ , }
\end{equation*}%
have the same distribution. Throughout this paper by ``isotropic'' we mean isotropic
in this sense.

We shall use the spherical coordinates
on $S^2$ i.e. $%
x =(\vartheta ,\varphi )$, where $0\leq \vartheta \leq \pi $ , $0\leq
\varphi <2\pi$. Also we shall use for $g\in SO(3)$ the
parameterization through the Euler angles $0\leq \alpha ,\gamma <2\pi $ and
$0\leq \beta \leq \pi $. Assuming the right-hand side exists, for each $\ell=1,2,\dots $we
can define the random vector
\begin{equation}
a_{\ell.}=\int_{S^{2}}T(x )Y_{\ell.}(x )\,dx   \label{basdef}
\end{equation}%
where $dx =\sin \vartheta\, d\varphi\, d\vartheta $ denotes the Lebesgue measure
on $S^{2}$ and $Y_{\ell.}$ the vector of spherical harmonics defined by
\begin{align}
Y_{\ell.}(\theta ,\varphi )& =\left( Y_{\ell\ell}(\theta ,\varphi ),\dots ,Y_{\ell,-\ell}(\theta ,\varphi
)\right) ^{\prime }\text{ ,}  \notag \\
Y_{\ell m}(\theta ,\varphi )& =\sqrt{\frac{2l+1}{4\pi }\frac{(\ell-m)!}{(\ell+m)!}}%
P_{\ell m}(\cos \theta )e^{im\varphi}\text{ , for }m\ge 0\text{ ,}  \notag \\
Y_{\ell m}(\theta ,\varphi )& =(-1)^{m}\overline{Y_{\ell,-m}(\theta ,\varphi )}\text{
, for }m<0\text{ };  \label{spharm}
\end{align}%
here 
$P_{\ell m}(\cos \theta )$
denotes the associated Legendre functions 
i.e.
\begin{align*}
P_{\ell m}(x)& =(-1)^{m}(1-x^{2})^{m/2}\frac{d^{m}}{dx^{m}}P_{\ell}(x)\text{ , }%
P_{\ell}(x)=\frac{1}{2^{\ell}l!}\frac{d^{\ell}}{dx^{\ell}}(x^{2}-1)^{\ell}, \\
m& =0,1,2,\dots ,\ell\text{ , }\ell=1,2,3,\dots .\text{ .}
\end{align*}%
A detailed discussion of the properties of the spherical harmonics can be
found in Varshalovich, Moskalev and Khersonskii \cite{MR1022665}, chapter 5. Our
purpose in this paper is to provide some characterizations of the
probability law of the vector $a_{\ell.}$ under the isotropy assumption. It is
a standard result that the functions $Y_{\ell m}$ $m=-\ell,\dots,\ell$ form a
basis for the vector space of functions on $S^2$ which are restrictions of homogeneous harmonic
polynomials of degree $\ell$. This vector space being invariant by the action of $SO(3)$,
for any $\ell$ and $g\in SO(3)$ there exist a $(2\ell+1)\times (2\ell+1)$ matrix $D^{\ell}(g)$ such that%
\begin{equation}
Y_{\ell.}(g x)=D^{\ell}(g)Y_{\ell}(x )\text{ ,}  \label{basrel}
\end{equation}%
This provides a representation of $SO(3)$ on ${\mathbb C}^{2\ell+1}$ which moreover is
irreducible (see Vilenkin and Klymik \cite{MR1220225}, \S 9.2.6 e.g.)
The matrices $D^{\ell}(g)$  are the so-called Wigner's
D-matrices whose entries are,
in terms of the Euler
angles
$$
D_{mm^{\prime }}^{\ell}(\alpha ,\beta ,\gamma ) =e^{-im\alpha }\,d_{mm^{\prime
}}^{\ell}(\beta )\,e^{-im^{\prime }\gamma },
$$
where
\begin{eqnarray*}
d_{mm^{\prime }}^{\ell}(\beta ) &=&(-1)^{\ell-m^{\prime} }\left[ (\ell+m)!\,(\ell-m)!\,(\ell+m^{%
\prime })\,!(\ell-m^{\prime })!\right] ^{1/2} \\
&\times& \sum_{k=0}^{\max (\ell-m,\ell-m^{\prime })}\frac{(\cos \frac{\beta }{2}%
)^{2k-m-m^\prime }(\sin \frac{\beta }{2})^{2\ell-2k-m-m^\prime }}{%
k!\,(\ell-m-k)!\,(\ell-m^{\prime }-k)!\,(m+m^{\prime }+k)!}\text{ .}
\end{eqnarray*}%
Note that $d_{mm^{\prime }}^{\ell}(0 )=\delta
_{m}^{m^\prime }$, where $\delta _{m}^{m\prime }$ denotes the
Kronecker delta function.
It is immediate that, the Lebesgue measure of $S^2$ being invariant by the action of
$SO(3)$,
\begin{equation}\label{invar1}
\begin{array}c
\displaystyle
a_{\ell.} =\int_{S^{2}}T(x )Y_{\ell.}(x )\,dx \overset{d}{=}%
\int_{S^{2}}T(gx )Y_{\ell.}(x )\,dx\cr
\displaystyle=\int_{S^{2}}T(x)D^{\ell}(g^{-1})Y_{\ell.}(x )\,dx
=D^{\ell}(g^{-1})a_{\ell.}\text{ ,}
\end{array}
\end{equation}
In particular, in coordinates, 
\begin{equation}
a_{\ell m}\overset{d}{=}\sum_{m^{\prime }=-\ell}^{\ell}a_{\ell m^{\prime }}D_{m^{\prime
}m}^{\ell}(g)\text{ .}  \label{invrot2}
\end{equation}%
It is well known that for an isotropic field which is continuous in mean square the
development (\ref{develop}) holds, the convergence being in $L^2$ (see \cite{MR1687092} e.g.).

It is also useful to point out that, because of (\ref{spharm}), the coefficients satisfy
the following identity
$$
a_{\ell m}=(-1)^m\overline{a_{\ell, -m}}.
$$
In particular $a_{\ell 0}$ is real.
\section{General properties of the spherical harmonics coefficients}
It is well-known that, for mean square continuous and isotropic random
fields, the spherical harmonics coefficients are orthogonal, i.e.
\begin{equation}
E[a_{\ell_{i}m_{i}}\overline{a_{\ell_{j}m_{j}}}]=\delta _{\ell_{i}}^{\ell_{j}}\delta
_{m_{j}}^{m_{i}}C_{\ell_{i}}\text{ , }  \label{orthocon}
\end{equation}
the sequence $\{ C_{\ell}\} _{\ell=1,2,\dots }$ denoting the angular power
spectrum of these fields.
For reasons of completeness we give now a proof of this fact. Actually our statement is
slightly stronger.
%
As usual we denote by $A^*$ the complex conjugate of the matrix $A$.
\begin{proposition}\label{prop1} \sl Assume $T$ isotropic. Then

a) for all $\ell$ such that $E[| a_{\ell.}| ^{2}]<\infty $,
\begin{equation*}
Ea_{\ell.}a_{\ell.}^{\ast }=C_{\ell}I_{2\ell+1}\text{ ,}
\end{equation*}%
where $I_{2\ell+1}$ denotes the $(2\ell+1)\times (2\ell+1)$ identity matrix

b) for all $\ell_{1},\ell_{2}$ such that $E[| a_{\ell_{1}.}|
^{2}]<\infty $, $E[| a_{\ell_{1}.}|
^{2}]<\infty $
\begin{equation*}
Ea_{\ell_{1}.}a_{\ell_{2}.}^{\ast }=0
\end{equation*}%
(in the sense of the $(2\ell_{1}+1)\times (2\ell_{2}+1)$
zero matrix).
\end{proposition}
{\it Proof } a) Let us denote by $\Gamma_\ell$ the covariance matrix of the random vector
$a_{\ell.}$. Since the vectors $a_{\ell.}$ and $D^{\ell}(g)a_{\ell.}$ have the same distribution,
they have the same covariance matrix. This gives
$$
\Gamma_\ell=D^{\ell}(g)\Gamma_\ell D^{\ell}(g)^*=D^{\ell}(g)\Gamma_\ell D^{\ell}(g)^{-1}
$$
Since $D^{\ell}$ is an irreducible representation of $SO(3)$, by Schur lemma
$\Gamma_\ell$
is of the form $C_{\ell}I_{2\ell+1}$.

b) The representations $D^{\ell_{1}}$ and $D^{\ell_{2}}$ are not equivalent for
$\ell_{1}\not=\ell_{2}$, having different dimensions. Therefore again by Schur lemma,
the identity
$$
Ea_{\ell_{1}.}a_{\ell_{2}.}^{\ast }=D^{\ell_{1}}(g)Ea_{\ell_{1}.}a_{\ell_{2}.}^{\ast }
D^{\ell_{1}}(g)^{-1}
$$
can hold only if the right hand side is the zero matrix.
\begin{remark} \rm We stress that (for strongly isotropic fields) Proposition
1 is strictly stronger than the standard result on mean square random
fields. Indeed, it is immediate to show that $ET^{2}<\infty $ implies $%
\sum_{l=1}^{\infty }E | a_{l.}| ^{2}<\infty ,$ on the other
hand, it is not difficult to find examples where mean square continuity
fails but the assumptions of Proposition 1 are fulfilled. Consider for
instance the field:%
\begin{equation}
T(x
)=\sum_{m=-\ell_{1}}^{\ell_{1}}a_{\ell_{1}m_{1}}Y_{\ell_{1}m_{1}}(x
)+\sum_{m=-\ell_{2}}^{\ell_{2}}a_{\ell_{2}m_{2}}Y_{\ell_{2}m_{2}}(x
)+\sum_{m=-\ell_{3}}^{\ell_{3}}b_{\ell_{3}m_{3}}Y_{\ell_{3}m_{3}}(x )\label{ex1}
\end{equation}
where $b_{\ell_{3}m_{3}}=\eta a_{\ell_{3}m_{3}}$; we assume that the $%
a_{\ell_{i}m_{i}}$'s $(i=1,2,3)$ satisfy (\ref{orthocon}) whereas $\eta $ is a
random variable with infinite variance (for instance a Cauchy). It is not
difficult to see that the field $T$ is properly defined
and strictly isotropic; although (\ref{ex1}) is clearly an artificial model,
some closely related field may be of interest for practical applications:
for instance in CMB data analysis it is often the case that the observed
field is a superposition of signal plus foreground contamination, and the
latter may be characterized by heavy tails at the highest multipoles (point
sources). In such cases, it is of an obvious statistical interest to know
that the standard properties of the spherical harmonics coefficients still
hold at least for the multipoles where foreground contamination is absent.
It is immediate to see that $ET^{2}=\infty ,$ whence the field cannot be
mean-square continuous; however (\ref{basdef}) is still properly defined for
$l=l_{1},l_{2}$ (simply exchange the integral with the finite sum), and
therefore Proposition 1 holds for these two vectors of spherical harmonics
coefficients.
\end{remark}
Our next proposition provides three further characterizations for the
spherical harmonics coefficients of isotropic fields.
\begin{proposition}\label{prop2} \sl Let $T$ be an isotropic
random field (not necessarily mean square continuous); then for all $a_{l.}$
such that $E | a_{\ell.} | ^{2}<\infty $,

a) for all $m=1,\dots ,\ell,$
\begin{equation*}
\func{Re} a_{\ell m}\overset{d}{=}\func{Im}a_{\ell m}\quad\mbox{ and }\quad \frac{\func{Re}a_{\ell m}%
}{\func{Im}a_{\ell m}}\sim Cauchy\text{  }
\end{equation*}

b) for all $m=1,2,\dots ,\ell,$ $\func{Re}a_{\ell m}$ and $\func{Im}a_{\ell m}$ are
uncorrelated, with variance $E(\func{Re}a_{\ell m})^{2}=E(\func{Im}%
a_{\ell m})^{2}=C_{\ell}/2.$

c) The marginal distribution of $\func{Re}a_{\ell m}$ , $\func{Im}a_{\ell m}$ is
always symmetric, that is,
\begin{equation*}
\func{Re}a_{\ell m}\overset{d}{=}-\func{Re}a_{\ell m}\text{ , }\func{Im}a_{\ell m}%
\overset{d}{=}-\func{Im}a_{\ell m}\text{ .}
\end{equation*}%
\end{proposition}
{\it Proof }
a) For $\beta =\gamma =0,$ (\ref{invrot2}) becomes%
\begin{equation*}
a_{\ell m}\overset{d}{=}e^{-im\alpha }a_{\ell m}
\end{equation*}%
for all $m=-\ell,\dots ,\ell$, $0\leq \alpha <2\pi$. This entails
\begin{equation}
\left(
\begin{array}{c}
\func{Re}a_{\ell m} \\
\func{Im}a_{\ell m}%
\end{array}%
\right) \overset{d}{=}\left(
\begin{array}{cc}
\cos \varphi  & \sin \varphi  \\
-\sin \varphi  & \cos \varphi
\end{array}%
\right) \left(
\begin{array}{c}
\func{Re}a_{\ell m} \\
\func{Im}a_{\ell m}%
\end{array}%
\right)
\label{rot1}
\end{equation}%
for all $m=-\ell,\dots ,\ell$, $0\leq \varphi <2\pi$.
Thus the vector $^t(\func{Re}a_{\ell m},\func{Im}a_{\ell m})$ has a distribution that is invariant
by rotations, that is, in polar coordinates, they can be written in of the form
\begin{equation}\label{rho1}
R \cos(\Theta)
\end{equation}
where $R$ is a random variable with values in $\mathbb{R}^+$, whereas $\Theta$ is uniform
in $[-\pi,\pi]$. This entails immediately that
$\arctan( \func{Re}a_{\ell m}/\func{Im}%
a_{\ell m}) \sim U(-\frac{\pi }{2},\frac{\pi }{2});$ the result follows
immediately.

b) This property is well-known if $T$ is mean square continuous. From (\ref{rho1}),
\begin{equation*}
E[\func{Re}a_{\ell m}\cdot \func{Im}a_{\ell m}]=\int_0^{+\infty}r^2d\mu_R(r)\int_{-\pi}^{\pi }\cos
\vartheta \sin \vartheta d\vartheta =0\text{ .}
\end{equation*}%

c) It suffices to take $\varphi _{0}=\pi $ in (\ref{rot1}).
\begin{remark}\rm It is interesting to note how Proposition \ref{prop1} implies that
no information can be derived on the statistical distribution of an
isotropic random field by the marginal distribution function of the ratios $%
( \func{Re}a_{\ell m}/\func{Im}a_{\ell m})$. On the other hand, it may
be possible to use these ratios to implement statistical tests of the
assumption of isotropy, an issue which has gained a remarkable empirical
relevance after the first release of the WMAP data in February 2003
\end{remark}
It is clear that if the field $T$ is gaussian, then the r.v.'s $(a_{\ell m})_{\ell,m}$
is a gaussian family. We prove now an independence result for this family of r.v.'s.
Thanks to Proposition \ref{prop1}, these
r.v.'s are uncorrelated, but one must be careful, since, in the case of complex r.v.'s,
 absence of correlation and
joint gaussian distribution does not imply independence.

\begin{proposition}\label{prop-if} \sl
For an isotropic gaussian random field the r.v.'s $a_{\ell m}$, $\ell=0,1,\dots$,
$m=0,\dots,\ell$ are independent.
\end{proposition}
{\it Proof} Let be $(\ell,m )\not=(\ell',m')$, $m>0$, $m'> 0$. Then $a_{\ell m}$
is uncorrelated with both $a_{\ell' m'}$ and $a_{\ell', -m'}=\overline{a_{\ell' m'}}$.
Thus
$$
E[a_{\ell m}\overline{ a_{\ell' m'}}]=0,\qquad E[a_{\ell m}a_{\ell' m'}]=E[a_{\ell m}
\overline{a_{\ell', -m'}}]=0
$$
and the statement follows from Lemma \ref{lemma}. If one at least among $m$ and $m'$ is equal to
$0$, then the r.v. $a_{\ell m}$ (or $a_{\ell' m'}$ is real and independence follows from
absence of correlation as for the real case.
\begin{lemma}\label{lemma}\sl
Let $Z_1$, $Z_2$ be complex r.v.'s, centered and jointly gaussian.
Then they are independent if and only if
\begin{equation}\label{complex1}
E[Z_1\overline Z_2]=0,\qquad E[Z_1 Z_2]=0
\end{equation}
\end{lemma}
{\it Proof} In one direction the statement is obvious. Let us assume that
(\ref{complex1}) are satisfied. Then, if we set $Z_k=X_k+iY_k$,
$k=1,2$, then
$$
\displaylines{
E[X_1X_2+Y_1Y_2]+i E[-X_1Y_2+Y_1X_2]=0\cr
E[X_1X_2-Y_1Y_2]+i E[X_1Y_2+Y_1X_2]=0\cr
}
$$
From these one obtains $E[X_1X_2]=0$, $E[Y_1Y_2]=0$, $E[X_1Y_2]=0$
and $E[Y_1X_2]=0$. This means that each of the r.v.'s $X_1,Y_1,X_2,Y_2$ is uncorrelated with
the other ones, so that, being jointly gaussian, they are independent.
\bigskip

\noindent Which is less obvious, is that the converse
also holds. The following is the main result of this paper.
%
\begin{theorem}\label{onlyif} \sl For an isotropic random field, let $\ell$ be such that $E | a_{\ell .}|
^{2}<\infty$. Then the coefficients $(a_{\ell 0},a_{\ell 1}\dots , a_{\ell \ell })$ are
independent if and only if they are gaussian.
\end{theorem}
{\it Proof } We just need to prove the ``only if'' part. Fix $m_1\ge 0$, $m_2\ge 0$, so
that the two complex r.v.'s $a_{\ell m_{1}}$ and $a_{\ell m_{2}}$ are independent.
Note that we are not assuming the independence of $\func{Re}a_{\ell m_{1}}$ and
$\func{Im}a_{\ell m_{1}}$ or of $\func{Re}a_{\ell m_{2}}$ and
$\func{Im}a_{\ell m_{2}}$. Thanks to (\ref{invar1}), the two vectors $a_{\ell .}$
and $D^\ell(g)a_{\ell .}$
have the same distribution. Thus the two r.v.'s
$$
L_1=\sum_{m'=-\ell}^\ell D^\ell_{m',m_1}(g)a_{\ell m'}\quad\mbox{and} \quad
L_2=\sum_{m'=-\ell}^\ell D^\ell_{m',m_2}(g)a_{\ell m'}
$$
having the same joint distribution as $a_{\ell m_1}$ and $a_{\ell m_2}$, are independent.
Fix $g$ so that the angles $\alpha $ such that $m\alpha
\neq k\pi $ for all integers $m,k$ and $\beta $ such that $%
d_{mm_{1}}^{\ell }(\beta )$ and $d_{mm_{2}}^{\ell }(\beta )$ are different from zero
for all $m=-\ell ,\dots ,\ell $ (note that such a $\beta $ certainly exists, because
the functions $d_{m_{2}m}^{\ell }(\beta )$, $m,m'=-\ell,\dots,\ell$ are analytic
and can vanish only at a finite number of values $\beta\in[0,\pi]$). For such a choice
of $g$,
thanks
to (\ref{invar1}),
\begin{equation*}
\left(
\begin{array}{c}
\func{Re}L_1 \\
\func{Im}L_1%
\end{array}%
\right) {=}\sum_{m=-\ell }^{\ell}d_{mm_{1}}^{\ell }(\beta )\left(
\begin{array}{cc}
\cos m_{1}\alpha  & \sin m_{1}\alpha  \\
-\sin m_{1}\alpha  & \cos m_{1}\alpha
\end{array}%
\right) \left(
\begin{array}{c}
\func{Re}a_{\ell m} \\
\func{Im}a_{\ell m}%
\end{array}%
\right)
\end{equation*}%
and%
\begin{equation*}
\left(
\begin{array}{c}
\func{Re}L_2 \\
\func{Im}L_1%
\end{array}%
\right) {=}\sum_{m=-\ell }^{\ell }d_{mm_{2}}^{\ell }(\beta )\left(
\begin{array}{cc}
\cos m_{2}\alpha  & \sin m_{2}\alpha  \\
-\sin m_{2}\alpha  & \cos m_{2}\alpha
\end{array}%
\right) \left(
\begin{array}{c}
\func{Re}a_{\ell m} \\
\func{Im}a_{\ell m}%
\end{array}%
\right)
\end{equation*}%
where the $2\times 2$ matrices on the right hand-sides are always full rank.
By the Skitovich-Darmois theorem below
(see Kagan, Rao and Linnik \cite{MR0346969} e.g.), it
follows that each of the vectors $(%
\func{Re}a_{\ell m},\func{Im}a_{\ell m})$ is bivariate Gaussian; as $(\func{Re}%
a_{\ell m},\func{Im}a_{\ell m})$ are uncorrelated and have the same variance by
Proposition \ref{prop1}, then $a_{\ell m}=\func{Re}a_{\ell m}+i\func{Im}a_{\ell m}$
is complex
Gaussian.
\begin{theorem} (Skitovich-Darmois) \sl Let $X_{1},\dots ,X_{r}$ be mutually
independent random vectors in $R^{n}.$ If the linear statistics
\begin{equation*}
L_{1}=\sum_{j=1}^{r}A_{j}X_{j},\qquad L_{2}=\sum_{j=1}^{r}B_{j}X_{j}\text{
,}
\end{equation*}%
are independent, for some real nonsingular $n\times n$ matrices $A_{j},B_{j},
$ $j=1,\dots ,r,$ then each of the vectors $X_{1},\dots ,X_{r}$ is normally
distributed.
\end{theorem}
In particular Theorem \ref{onlyif} implies that if an isotropic random field is mean square
continuous and the coefficients $a_{\ell m}, \ell=0,1,\dots, m=0,\dots \ell$ are independent, then
it is gaussian.
\begin{remark}\rm Proposition \ref{onlyif} shows that it is not possible to generate
isotropic random fields by sampling non-Gaussian, independent complex-valued
random variables $a_{\ell m},$ $m=-\ell,\dots,\ell$. This fact shows that, apart from the
gaussian case, it is not easy to sample a random field by simulating the values of
the random coefficients $a_{\ell m}$. In particular sampling independent values of
the $a_{\ell m}$'s with distributions other than the gaussian gives not rise to an isotropic
random field.

We wish also to point out that it is indeed possible to construct a non-gaussian
random field by choosing the random
coefficients $a_{\ell m}$, $\ell=0,1,\dots$ $m=0,\dots,\ell$
independent and with an arbitrary distribution. If they satisfy
the conditions of Proposition \ref{prop1} and the series (\ref{develop}) converges, they
certainly define a random field on $S^2$. But Theorem \ref{onlyif} states that such a
field
cannot be isotropic. In particular Theorem \ref{onlyif} does not follow by any means from
the
central limit theorem.
\end{remark}

\noindent{\bf Acknowledgment.} The authors wish to thank Professor R. Varadarajan for a
very illuminating discussion and useful remarks.
\bibliography{bibbase}

\begin{thebibliography}{10}

\bibitem{AGMNW:04}
{\sc M.~Arjunwadkar, C.~R. Genovese, C.~J. Miller, R.~C. Nichol, and
  L.~Wasserman}, {\em Nonparametric inference for the cosmic microwave
  background.}, Statist. Sci., 19 (2004), pp.~308--321.

\bibitem{MR2106384}
{\sc N.~Bartolo, E.~Komatsu, S.~Matarrese, and A.~Riotto}, {\em
  Non-{G}aussianity from inflation: theory and observations}, Phys. Rep., 402
  (2004), pp.~103--266.

\bibitem{CoMa:2001}
{\sc C.~R. Contaldi and J.~Magueijo}, {\em Generating non-gaussian maps with a
  given power spectrum and bispectrum}, Phys. Rev. D, 63 (2001).

\bibitem{cruz:05}
{\sc M.~Cruz, E.~Marinez-Gonzalez, P.~Vielva, and L.~Cayon}, {\em Detection of
  a non-gaussian spot in {WMAP}}, Monthly Notices of R. Astronomical Society,
  356 (2005), pp.~29--40.

\bibitem{DCB:2003}
{\sc O.~Dor\'e, S.~Colombi, and F.~Bouchet}, {\em Probing {CMB} non
  {Gaussianity} using local curvature}, Monthly Notices of R. Astronomical
  Society, 344 (2003), pp.~905--916.

\bibitem{JSDAF:05}
{\sc J.~Jin, J.-L. Starck, D.~Donoho, N.~Aghanim, and O.~Forni}, {\em
  Cosmological non-gaussian signature detection: Comparing performance of
  different statistical tests.}, Eurasip J. Appl. Signal Processing,  (2005),
  p.~forthcoming.

\bibitem{MR0346969}
{\sc A.~M. Kagan, Y.~V. Linnik, and C.~R. Rao}, {\em Characterization problems
  in mathematical statistics}, John Wiley \& Sons, New York-London-Sydney,
  1973.
\newblock Translated from the Russian by B. Ramachandran, Wiley Series in
  Probability and Mathematical Statistics.

\bibitem{MR1687092}
{\sc N.~Leonenko}, {\em Limit theorems for random fields with singular
  spectrum}, vol.~465 of Mathematics and its Applications, Kluwer Academic
  Publishers, Dordrecht, 1999.

\bibitem{MA:2005}
{\sc D.~Marinucci}, {\em Testing for non-gaussianity on cosmic microwave
  background radiation: A review}, Statist. Sci., 19 (2004), p.~294–307.

\bibitem{MR2065205}
{\sc D.~Marinucci and M.~Piccioni}, {\em The empirical process on {G}aussian
  spherical harmonics}, Ann. Statist., 32 (2004), pp.~1261--1288.

\bibitem{Pa:2004}
{\sc C.~Park}, {\em Non gaussian signatures in the temperature fluctuations
  observed by the {Wilkinson Microwave Anisotropy Probe}}, Monthly Notices of
  R. Astronomical Society, 349 (2004), pp.~313--320.

\bibitem{MR1022665}
{\sc D.~A. Varshalovich, A.~N. Moskalev, and V.~K. Khersonski{\u\i}}, {\em
  Quantum theory of angular momentum}, World Scientific Publishing Co. Inc.,
  Teaneck, NJ, 1988.
\newblock Irreducible tensors, spherical harmonics, vector coupling
  coefficients, $3nj$ symbols, Translated from the Russian.

\bibitem{MR1220225}
{\sc N.~J. Vilenkin and A.~U. Klimyk}, {\em Representation of {L}ie groups and
  special functions. {V}ol. 2}, vol.~74 of Mathematics and its Applications
  (Soviet Series), Kluwer Academic Publishers Group, Dordrecht, 1993.
\newblock Class I representations, special functions, and integral transforms,
  Translated from the Russian by V. A. Groza and A. A. Groza.

\end{thebibliography}
\bibliographystyle{siam}
\end{document}